\documentclass{amsart}
\title{Short proofs for $q$-Raabe formula and integrals for Jacobi theta functions}
\usepackage{amssymb,amsmath,amsthm,epsfig,graphics,latexsym}
\usepackage{enumerate}
 \theoremstyle{definition}

  \theoremstyle{plain}

  \newtheorem{theorem}    {Theorem}

  \DeclareMathOperator{\li}{Li}
\begin{document}
  \author{Mohamed El Bachraoui}
 \email{melbachraoui@uaeu.ac.ae}
  \begin{abstract}
  We shall answer a question of 
  Mez\H{o} on the $q$-analogue of the Raabe's integral formula for $0<q<1$ and we shall  evaluate an integral involving the first theta function. Moreover, we will reproduce short proofs for some identities of  
  Mez\H{o}.
  \end{abstract}
  \date{\textit{\today}}
  \maketitle
\section{Introduction}
Recall that for a complex number $q$ and a complex variable $a$, the $q$-shifted factorials are given by
\[
(a;q)_0= 1,\quad (a;q)_n = \prod_{i=0}^{n-1}(1-a q^i),\quad
(a;q)_{\infty} = \lim_{n\to\infty}(a;q)_n =\prod_{i=0}^{\infty}(1-a q^i)\quad (|q|<1)
\]
and recall the dilogarithm function
\[
\li(z) = \sum_{n=1}^{\infty} \frac{z^n}{n^2}
\]
which for $z=1$ evaluates to $\zeta(2)$.
There are known two $q$-analogues of the gamma function which both were first
introduced by Jackson~in~\cite{Jackson}. The first one is:
\begin{equation}\label{q-gamma-1}
\Gamma_q (x) = \frac{(q;q)_{\infty}}{(q^x;q)_{\infty}} (1-q)^{1-x}\quad (0<q<1)
\end{equation}
and the second one is:
\begin{equation}\label{q-gamma-2}
\Gamma_q (x) = \frac{(q^{-1};q^{-1})_{\infty}}{(q^{-x};q^{-1})_{\infty}} (q-1)^{1-x}q^{{x\choose 2}}
\quad (q>1).
\end{equation}
For more details on the version of $\Gamma_q(x)$ for $0<q<1$ we refer to \cite{Askey-1, Askey-2}
and on the version of $\Gamma_q(x)$ for $q>1$ we refer to \cite{Moak}.
Raabe~\cite{Raabe} gave the following integral
\begin{equation}\label{Raabe-General}
\int_{0}^1 \log \Gamma(x+t) \,dx = \log\sqrt{2\pi} + t\log t - t\qquad (t\geq 0)
\end{equation}
which as $t\to 0+$ implies
\begin{equation}\label{Raabe-Special}
\int_{0}^1 \log \Gamma(x) \,dx = \log\sqrt{2\pi}.
\end{equation}
Recently, Mez\H{o} found $q$-analogues for both (\ref{Raabe-General}) and (\ref{Raabe-Special}) as in the following theorem.
\begin{theorem} \emph{(Mez\H{o}~\cite[Theorem 2]{Mezo})} \label{Mezo-Thm-2}
If $q>1$, then for any $t>0$,
\begin{equation}\label{q-Raabe-General}
\begin{split}
\int_{0}^1 \log \Gamma_q(x+t) \,dx &=
\log C_q - \frac{1}{2 q^t \log q} \Big( \frac{1-q^t}{1-q^{-t}} \big( 2\li_2(q^{-t}) + \log^2(1-q^{-t}) \big) \\
& + 2\frac{1-q^t}{1-q^{-t}} \log\frac{1-q}{1-q^{t}} \log(1-q^{-t}) - q^t \log^2 \frac{1-q}{1-q^{t}} \Big),
\end{split}
\end{equation}
where
\[
C_q = q^{-\frac{1}{12}} (q-1)^{ \frac{1}{2}- \frac{\log(q-1)}{2\log q} } (q^{-1};q^{-1})_{\infty}.
\]
In particular, if $t\to 0$, then
\begin{equation}\label{q-Raabe-Special}
\int_{0}^1 \log \Gamma_q(x) \,dx = \frac{\zeta(2)}{\log q} + \log\sqrt{ \frac{q-1}{\sqrt[6]{q}} } + \log(q^{-1};q^{-1})_{\infty}.
\end{equation}
\end{theorem}
\noindent
To find these formulas, Mez\H{o} needed to evaluate the integral $\int_{0}^1 \zeta_q(s,x+t) \,dx$ of the 
$q$-Hurwitz zeta function and made an appeal to a result by
 Kurokawa~and~Wakayama~\cite{Kurokawa-Wakayama}. 
 However, these ideas seem not to apply directly to the
 case $0<q<1$ and  therefore the author asked how identities 
 (\ref{q-Raabe-General}) and (\ref{q-Raabe-Special}) look like in the latter case.
 In this note we will answer Mez\H{o}'s question as in the following theorem.
 \begin{theorem} \label{Own q-Raabe}
If $0<q<1$ and $t\geq 0$, then
\[
\int_{0}^1 \Gamma_q(x+t) \,dx = \big(\frac{1}{2}-t \big) \log (1-q) - \frac{1}{\log q} \li_2(q^t) + \log(q;q)_{\infty}.
\]
In particular, if $t=0$, then
\[
\int_{0}^1 \Gamma_q(x) \,dx = \frac{1}{2}\log (1-q) - \frac{\zeta(2)}{\log q} + \log(q;q)_{\infty}.
\]
\end{theorem}
 Next, using the same approach we will reproduce a short, elementary proof for 
 Mez\H{o}'s Theorem~\ref{Mezo-Thm-2}.
 In fact, Mez\H{o}'s main result in \cite{Mezo} is the following theorem
 involving the Jacobi's fourth theta function
 \[
 \theta_4 (x,q) = \sum_{n=-\infty}^{\infty} (-1)^n q^{n^2} e^{2nix}.
 \]
 \begin{theorem}\emph{(Mez\H{o}~\cite[Theorem 1]{Mezo})} \label{Mezo-Thm-1}
If $0<q<1$ is real, then
\[
\int_{-\frac{1}{2}\log q}^{\frac{1}{2}\log q} \log \theta_4(ix,q) \,dx = \zeta(2)+\log q \cdot \log(q^2;q^2)_{\infty}.
\]
\end{theorem}
To prove the previous result, the author among other things made an appeal to 
Theorem~\ref{Mezo-Thm-2}. In this paper we will provide a short proof for Theorem~\ref{Mezo-Thm-1}.
Furthermore, we will prove the following related theorem on the Jacobi's first theta function
\[
\theta_1 (x,q) = \sum_{n=-\infty}^{\infty} (-1)^{n-1/2} q^{(n+1/2)^2} e^{(2n+1)ix} .
\]
\begin{theorem}\label{int-theta1}
If $0<q<1$ is real, then
\[
\int_{0}^{\log q} \log \theta_1(ix,q) \,dx = \zeta(2) + \log q\cdot \log(q^2;q^2)_{\infty}.
\]
\end{theorem}
\section{Proof of Theorem \ref{Own q-Raabe}}
We start by the second identity.
It is clear that
\begin{equation}\label{log}
\log (q^x;q)_{\infty} = \log\prod_{k=0}^{\infty}(1-q^{x+k}) = \sum_{k=0}^{\infty} \log(1-q^{x+k}) 
\end{equation}
and that for each $k=0,1,\ldots$
\begin{equation}\label{indef-int-1}
\int \log(1-q^{x+k}) \,dx = \int \sum_{n=1}^{\infty}\frac{(q^{x+k})^n}{n}= -\frac{1}{\log q} \li_2 (q^{x+k}) + Constant.
\end{equation}
Then for each $k=0,1,\ldots$
\[
\int_{0}^1 \log(1-q^{x+k}) \,dx = \frac{1}{\log q} \sum_{n=1}^{\infty} \frac{q^{kn}}{n^2}- \frac{1}{\log q}\sum_{n=1}^{\infty} \frac{q^{(k+1)n}}{n^2}
\]
which combined with (\ref{log}) yields
\begin{equation}
\int_{0}^1 \log(q^x;q)_{\infty} \,dx = \sum_{k=0}^{\infty}\int_{0}^1 \log(1-q^{x+k}) \,dx = \frac{\zeta(2)}{\log q}.
\end{equation}
Now using
 the previous integral and the definition  (\ref{q-gamma-1}) we find
\[
\begin{split}
\int_{0}^1 \log\Gamma_q(x) \,dx  &= \int_{0}^1 \Big( \log (q;q)_{\infty} + (1-x) \log(1-q) - \log(q^x;q)_{\infty} \Big) \,dx \\
&= \log (q;q)_{\infty} + \frac{1}{2} \log(1-q) - \frac{\zeta(2)}{\log q},
\end{split}
\]
as desired.
As to the first identity,
the substitution rule applied to the indefinite integral (\ref{indef-int-1}) gives
\[
\begin{split}
\int_{0}^{1} \log(1-q^{x+t+k}) \,dx &= \int_{t}^{t+1} \log(1-q^{u}) \,du \\
&= \frac{1}{\log q} \sum_{n=1}^{\infty} \frac{q^{(t+k)n}}{n^2}- \frac{1}{\log q}\sum_{n=1}^{\infty} \frac{q^{(t+k+1)n}}{n^2},
\end{split}
\]
from which we get
\[
\int_{t}^{t+1} \log(q^{u};q)_{\infty} \,du = \frac{1}{\log q} \li_2(q^t).
\]
Now by definition~(\ref{q-gamma-1}) and the previous integral we conclude that
\[
\begin{split}
\int_{0}^{1} \log\Gamma_q(x+t) \,dx &= \int_{t}^{t+1} \log\Gamma_q(u) \,du \\
&= \int_{t}^{t+1} \Big( \log(q;q)_{\infty} + (1-u)\log (1-q) - \log(q^u;q)_{\infty} \Big) \,du \\
&= \log(q;q)_{\infty} + (\frac{1}{2}-t) \log(1-q) - \frac{1}{\log q} \li_2 (q^t).
\end{split}
\]
This completes the proof.
\section{A short proof for Theorem \ref{Mezo-Thm-2}}
It is easy to check that if $q>1$, then
\begin{equation}\label{key-identity}
\Gamma_q(x) = \Gamma_{q^{-1}}(x) q^{{x-1\choose 2}}.
\end{equation}
Thus with the help of Theorem~\ref{Own q-Raabe} we have
\[
\begin{split}
\int_0^1 \Gamma_q(x) \,dx &= \int_0^1 \Big( \log \Gamma_{q^{-1}}(x) + \frac{(x-1)(x-2)}{2}\log q \Big)
\,dx \\
&= \frac{\zeta(2)}{\log q} -\frac{1}{12}\log q+\frac{1}{2}\log(q-1) + \log (q^{-1};q^{-1})_{\infty} \\
&= \frac{\zeta(2)}{\log q} + \log\sqrt{\frac{q-1}{\sqrt[6]{q}}} + \log (q^{-1};q^{-1})_{\infty}.
\end{split}
\]
This proves identity~(\ref{q-Raabe-Special}). As to idenity~(\ref{q-Raabe-General}) we similarly get
\[
\begin{split}
\int_0^1 \Gamma_q(x+t) \,dx &= \int_0^1 \Big( \log \Gamma_{q^{-1}}(x+t) +
\frac{(x+t-1)(x+t-2)}{2}\log q \Big) \,dx \\
&= \log (q^{-1};q^{-1})_{\infty} + \big(\frac{1}{2}-t\big)\log(1-q^{-1}) - \frac{\li_2(q^{-t})}{\log q^{-1}} \\
& \quad\quad +\big( \frac{5}{12}+ \frac{t(t-2)}{2} \big) \log q \\
&= \log (q^{-1};q^{-1})_{\infty} + \frac{\li_2(q^{-t})}{\log q} + \big(\frac{1}{2}-t\big)\log(q-1) \\
&\quad\quad - \big(\frac{1}{2}-t\big) \log q + \big( \frac{5}{12}+ \frac{t(t-2)}{2} \big) \log q \\
&= \log (q^{-1};q^{-1})_{\infty} + \frac{\li_2(q^{-t})}{\log q} + \big(\frac{1}{2}-t\big)\log(q-1) \\
&\quad\quad +\big( \frac{t^2}{2}-\frac{1}{12} \big) \log q,
\end{split}
\]
which by a straightforward but long calculation can be verified to agree with the right-hand-side of the
formula~(\ref{q-Raabe-General}).
\section{A short proof for Theorem \ref{Mezo-Thm-1}}
By the triple product identity,  \cite{Gasper-Rahman, Wittaker-Watson}, we have
\[
\theta_4 (ix,q) = (q e^{-2x};q^2)_{\infty} (q e^{2x};q^2)_{\infty} (q^2;q^2)_{\infty}.
\]
Then
\begin{equation}\label{help-int-theta4}
\int_{-\frac{1}{2}\log q}^{\frac{1}{2}\log q} \log q \theta_4(ix,q) \,dx = 
\int_{-\frac{1}{2}\log q}^{\frac{1}{2}\log q} \Big( \log(q e^{-2x};q^2)_{\infty} + 
\log(q e^{2x};q^2)_{\infty} + \log(q^2;q^2)_{\infty} \Big) \,dx.
\end{equation}
It is easy to check that
\[
\begin{split}
\int_{-\frac{1}{2}\log q}^{\frac{1}{2}\log q} \log(1-q^{2k+1}e^{-2x}) \,dx &= 
- \sum_{n=1}^{\infty}\frac{(q^{2k+1})^n}{n} \int_{-\frac{1}{2}\log q}^{\frac{1}{2}\log q} e^{-2nx} \,dx \\
&= \frac{1}{2}\sum_{n=1}^{\infty}\frac{(q^{2k})^n}{n^2} - \frac{1}{2}\sum_{n=1}^{\infty}\frac{(q^{2k+2})^n}{n^2}
\end{split}
\]
implying that 
\begin{equation}\label{help-int-theta4-1}
\int_{-\frac{1}{2}\log q}^{\frac{1}{2}\log q} \log(q e^{-2x};q^2)_{\infty} = \frac{\zeta(2)}{2}.
\end{equation}
A similar argument shows that
\begin{equation}\label{help-int-theta4-2}
\int_{-\frac{1}{2}\log q}^{\frac{1}{2}\log q} \log(q e^{2x};q^2)_{\infty} \,dx = \frac{\zeta(2)}{2}.
\end{equation}
Now putting (\ref{help-int-theta4-1}) and (\ref{help-int-theta4-2}) in (\ref{help-int-theta4}) gives the desired integral.
\section{Proof of Theorem \ref{int-theta1}}
By the triple product identity (see \cite{Gasper-Rahman, Wittaker-Watson})
\[
\theta_1(ix,q) = (e^{-2x} q^2;q^2)_{\infty} (e^{2x} ;q^2)_{\infty} (q^2;q^2)_{\infty} 
\]
and therefore,
\begin{equation}\label{help-int-theta1}
\int_{0}^{\log q} \log  \theta_1(ix,q) \,dx =
\int_{0}^{\log q} \Big( \log(e^{-2x}q^2;q^2)_{\infty} +
\log(e^{2x};q^2)_{\infty} + \log(q^2;q^2)_{\infty} \Big) \,dx.
\end{equation}
Following the same ideas of our proof for Theorem~\ref{Mezo-Thm-1} above, we get
\[
\int_{0}^{\log} (e^{-2x} q^2;q^2)_{\infty} \,dx = \int_{0}^{\log q} (e^{2x};q^2)_{\infty} \,dx=
\frac{\zeta(2)}{2}.
\]
Now putting together in (\ref{help-int-theta1}) gives the desired identity.

%

\begin{thebibliography}{9}
\bibitem{Andrews-1}
G. E. Andrews,
\emph{The theory of partitions, Vol. 2},
Cambridge University Press, 1984.
%
\bibitem{Askey-1}
R. Askey,
\emph{The $q$-gamma and $q$-beta functions},
Appl. Anal. 8 (2) (1978/1979), 125-141.
%
\bibitem{Askey-2}
R. Askey,
\emph{Ramanujan's extensions of the gamma and beta functions},
Amer. Math. Monthly 87 (5) (1980), 346-359.
%
\bibitem{Gasper-Rahman}
G. Gasper and M. Rahman,
\emph{Basic Hypergeometric series},
Cambridge University Press, 2004.
%
\bibitem{Jackson}
F. H. Jackson,
\emph{The basic gamma-function and the elliptic functions},
Proc. Roy. Soc. London Ser. A. 76 (1905), 127-144.
%
\bibitem{Kurokawa-Wakayama}
N. Kurokawa and N. Wakayama,
\emph{Generalized zeta regularizations, quantum class number formulas, and Appell's $\mathcal{O}$-functions},
Ramanujan J. 10 (2005), 291-303.
%
\bibitem{Mezo}
I. Mez\H{o},
\emph{A $q$-Raabe formula and an integral of the fourth Jacobi theta function},
J. Number Theory 133 (2013, 692-704.
%
\bibitem{Moak}
D. S. Moak,
\emph{The $q$-gamma function for $q>1$},
Aequationes Math. 20 (1980), 278-285.
%
\bibitem{Raabe}
J. L. Raabe,
\emph{Angen\"{a}herte Bestimmung der Factorenfolge $1\cdot 1\cdot 2\cdot 3 \cdot 4 \cdot 5 \cdots n = \Gamma(1+n)=\int x^n e^{-x} dx$, wenn $n$ eine sehr grosse Zahl ist},
J. Reine Angew. Math. 25 (1840), 146-159.
%
\bibitem{Wittaker-Watson}
E.T. Wittaker and G.N. Watson,
\emph{A course of modern analysis},
Cambridge University Press, 1996.

%
%
\end{thebibliography}
\end{document}